\documentclass[12pt]{amsart}
\usepackage{a4}
\usepackage[cp1250]{inputenc}
\usepackage{amssymb}

\usepackage[T1]{fontenc}

\usepackage{graphicx}

\usepackage{color}

\newcounter{liste}

\newenvironment{lister}{\begin{list}%
{{\rm(\alph{liste})\hfill}}{%
\usecounter{liste}%
\setlength{\labelsep}{0ex}%
\setlength{\labelwidth}{\parindent}%
\setlength{\leftmargin}{\parindent}%
\addtolength{\leftmargin}{\parindent}%
\addtolength{\labelwidth}{\parindent}%
\addtolength{\topsep}{1.5\topsep}%
}}{\end{list}}

\newenvironment{lister*}{\begin{list}%
{{\hfill$\bullet$\hfill}}{%
\setlength{\labelsep}{0ex}%
\setlength{\labelwidth}{\parindent}%
\setlength{\leftmargin}{\parindent}%
\addtolength{\topsep}{1.5\topsep}%
}}{\end{list}}

\newenvironment{plainlist*}{\begin{list}%
{}{\setlength{\labelsep}{0ex}%
\setlength{\labelwidth}{0ex}%
\setlength{\leftmargin}{0ex}%
\addtolength{\topsep}{1.5\topsep}%
}}{\end{list}}

\newtheorem{theo}{Theorem}
\newtheorem*{cor*}{Corollary}
\newtheorem{prop}{Proposition}
\newtheorem{prob}{Problem}
\newtheorem*{rems*}{Remarks}

\newcommand{\id}{\mathop{\mathrm{id}}\nolimits}
\newcommand{\diam}{\mathop{\mathrm{diam}}\nolimits}

\newcommand{\cl}{\mathop{\mathrm{cl}}\nolimits}
\renewcommand{\int}{\mathop{\mathrm{int}}\nolimits}

\title[End points in chainable continua]{Feeding and killing end points in chainable continua}

\author[J.\ Krzempek]{Jerzy Krzempek}

\address{Faculty of Applied Mathematics\\Silesian University of Technology\\Kaszubska 23\\44-100 Gliwice\\Poland}

\email{jkrzempek.math@gmail.com}



\dedicatory{Dedicated to Professor Jerzy Mioduszewski, my amazing mathematical grandfather}

\begin{document}

\begin{abstract}
Using the classical technique of condensation of singularities, we prove that, for every zero-dimensional, complete separable metric
space $G$, there exists a Suslinian, chainable metric continuum whose set of end points is homeomorphic to~$G$. This answers a
question posed by R.\ Adikari and W.\ Lewis in [Houston J.\ Math.\ 45 (2019), no.\ 2, pp.\ 609--624].
\end{abstract}

\maketitle

{\em All spaces considered in this paper are metric and separable, and all maps are continuous.} The terminology follows R.\ Engelking
\cite{eng-89,eng-95} and S.\ B.\ Nadler \cite{nad-92}.

A continuum (= non-empty, connected compact space) $X$ is said to be {\bf chainable}\/ if, for every $\varepsilon>0$, $X$ has a finite
covering by open sets $U_1,\ldots,U_n$ such that $\diam U_i<\varepsilon$ and $U_i\cap U_j\neq\emptyset$ iff $|i-j|\leq 1$ for
$i,j=1,\ldots,n$. An element $p\in X$ is called an {\bf end point}\/ of $X$ if we can moreover assume that $p\in U_1$.

R.\ H.\ Bing initiated investigation of end points of chainable continua in the seminal 1951 paper \cite{bing-51}. In the recent
detailed study \cite{adl-19} R.\ Adikari and W.\ Lewis remark {\em inter alia}\/ that the set of end points of a non-degenerate,
hereditarily decomposable\footnote{A continuum is {\bf decomposable}\/ if it is the union of two proper subcontinua, otherwise it is
{\bf indecomposable}. A continuum is {\bf hereditarily decomposable}\/ if each of its non-degenerate subcontinua is decomposable.}
chainable continuum is non-empty, nowhere dense, $G_\delta$ in the continuum, and does not contain any non-degenerate connected
subset. On the other hand, given any zero-dimensional compact space~$F$, they use inverse limits and construct an hereditarily
decomposable chainable continuum whose set of end points is homeo\-morph\-ic to~$F$. In the present note we follow this trail and
answer their Questions 4 and~5 by proving

\begin{theo}\label{main}
For every zero-dimensional complete space~$G$, there exists a Suslin\-ian\footnote{ A continuum is said to be {\bf Suslinian} if every
family of its pairwise disjoint non-degenerate subcontinua is countable. Since any non-degenerate indecomposable continuum has
$2^{\aleph_0}$ pairwise disjoint composants, {\em every Suslinian continuum is hereditarily decomposable.}} chainable continuum whose
set of end points is homeomorphic to~$G$.
\end{theo}

\begin{cor*}
There is a Suslinian chainable continuum whose set of end points is homeo\-morph\-ic to the set of irrationals.\qed
\end{cor*}

\section{Preliminaries: Atomic maps, chainable continua, and condensation of singularities}

The following characterisation elucidates the notion of an end point.

\begin{theo}[Bing {\cite[Theorem 13]{bing-51}}]\label{bi}
For a point $p$ in a chainable continuum $X$, the following statements are equivalent.
\begin{lister}
\item $p$ is an end point of $X$.
\item For every pair of subcontinua $P,Q\ni p$ of $X$, either $P\subset Q$ or $Q\subset P$.\qed
\end{lister}
\end{theo}

Let us adopt the natural generalisation: If a point $p$ of a continuum $X$, not necessarily chainable, satisfies the statement (b) of
Bing's Theorem 2, then we shall call $p$ an {\bf end point} of~$X$.

A continuum $X$ is said to be {\bf unicoherent} if, for each pair of its subcontinua $P$,~$Q$ such that $X=P\cup Q$, the intersection
$P\cap Q$ is connected. $X$ is {\bf hereditarily unicoherent} if each of its subcontinua is unicoherent. $X$ is called a {\bf triod}
if it contains a subcontinuum $P$ such that $X\setminus P$ if the union of three non-empty sets each two of which are mutually
separated in~$X$. A continuum is {\bf atriodic} if it does not contain any triod. A subcontinuum $P$ of a space $X$ is {\bf terminal}
in~$X$ if, for every subcontinuum $Q\subset X$ that meets $P$, either $P\subset Q$ or $Q\subset P$. A map $f\colon X\to Y$ is said to
be {\bf monotone} [respectively: {\bf atomic}] if each of its point-inverses $f^{-1}(y)$, $y\in Y$, is a subcontinuum [respectively:
terminal subcontinuum] of~$X$ (hence, $f$ is a map onto~$Y$ as any continuum is non-empty). The following well-known proposition is
easily proved.

\begin{prop}
Suppose that $f\colon X\to Y$ is an atomic map, and $P$ is a subcontinuum of $X$ with $f(P)$ non-degenerate. Then
\begin{lister}
\item The formula $P=f^{-1}(f(P))$ holds.
\item If $f(P)$ is a decomposable continuum {\em [}respectively: indecomposable continuum, triod, unicoherent continuum\/{\em ]}, then so is~$P$.
\item If $f(P)$ and each of the point-inverses $f^{-1}(y)$, $y\in f(P)$, is an atriodic {\em [}respectively: hereditarily decomposable,
hereditarily unicoherent\/{\em ]} continuum, then so is~$X$.
\item If $f(P)$ is a Suslinian continuum and there are at most countably many non-degenerate point-inverses $f^{-1}(y)$, $y\in f(P)$,
which are moreover Suslinian, then $P$ itself is Suslinian.
\item An element $p\in P$ is an end point of $P$ if, and only if $p$ is an end point of the continuum $f^{-1}(f(p))$ and the image $f(p)$ is an end point of $f(P)$.\qed
\end{lister}
\end{prop}

We shall use the following characterisation of chainable continua.

\begin{theo}[J.\ B.\ Fugate\footnote{In \cite[p.\ 384]{fug-69} Fugate has a notion of a terminal subcontinuum, which is different from ours.
We would prefer to name his terminal subcontinuum an end subcontinuum because its definition is rather a generalisation of the
statement (b) of Theorem \ref{bi}.} {\cite[Theorem 2]{fug-69}}] A continuum $X$ is chainable if, and only if $X$~is atriodic,
hereditarily unicoherent, and each indecomposable subcontinuum of~$X$ is chainable.\qed
\end{theo}

\begin{cor*}[T.\ Maćkowiak {\cite[Proposition 11(iii)]{mac-85}}]
Suppose that $f\colon X\to Y$ is an atomic map, and $P$ is a subcontinuum of $X$ with $f(P)$ non-degenerate. Then
\begin{lister}\setcounter{liste}{5}
\item If $f(P)$ is chainable and hereditarily decomposable, and each of the point-inverses $f^{-1}(y)$, $y\in f(P)$, is chainable,
then $P$~is chainable.
\end{lister}
\end{cor*}

\begin{proof}
$P=f^{-1}(f(P))$ by Proposition 1(a). If $f(P)$ is hereditarily decomposable, then by Proposition 1(a, b) every indecomposable
subcontinuum of $P$ is contained in a point-inverse under~$f$. Since {\em every chainable continuum is atriodic and hereditarily
unicoherent} (see \cite[Theorems 12.2 and 12.4]{nad-92}), Proposition 1(c) and Fugate's theorem imply~(f).
\end{proof}

Condensation of singularities has its origin in Z.\ Janiszewski's note \cite{jan-13}, where he constructed a non-degenerate continuum
that contains no arc. The technique was developed by G.\ T.\ Whyburn \cite{why-30}, R.\ D.\ Anderson and G.\ Choquet \cite{ach-59},
who put inverse limits to work, and was further augmented by Ma\'{c}kowiak \cite{mac-85} and E.\ Pol and M.\ Re\'{n}ska \cite{por-03}.

Janiszewski's arcless continuum, L.\ E.\ J.\ Brouwer's, Janiszewski's, and K.\ Yo\-ne\-ya\-ma's indecomposable continua (let us think
of the bucket-handle continuum), or B.\ Kna\-s\-ter's hereditarily indecomposable continuum (called by E.\ E.\ Moise the pseudo-arc)
were only possible in Cantor's paradise. It is interesting to see how these examples were simultaneously represented as 
inverse limits: Anderson and Choquet, the tent map for the bucket-handle, J.\ R.\ Isbell's proof \cite{isb-59} that {\em every
chainable continuum is an inverse limit with bonding maps from $[0,1]$ to itself,}\/ and J.\ Mio\-du\-szew\-s\-ki's functional
approach to chainable continua and the pseudo-arc \cite{mio-62,mio-64}. All the mentioned papers by Anderson, Choquet, Isbell, and
Mioduszewski were published between 1959 and 1964, and they displayed categorical thinking.

\begin{theo}[Pol and Re\'{n}ska {\cite[Theorem 3.2]{por-03}}]
Suppose that $X$~is a continuum, $(A_i)_{i=1}^\infty$ is a sequence of pairwise disjoint, at most zero-dimensional closed subsets of\/
$X$, and\/ $(Z_i)_{i=1}^\infty$ is any sequence of compact spaces. If each $Z_i$ admits a monotone map onto $A_i$, then there exists a
continuum $L(X,Z_i,A_i)$ with a surjective map $pr\colon L(X,Z_i,A_i)\to X$ which have the following properties.
\begin{lister}
\item $pr$ is atomic.
\item The restriction $pr|pr^{-1}(X\setminus\bigcup_{i=1}^\infty
A_i)\colon pr^{-1}(X\setminus\bigcup_{i=1}^\infty A_i)\to X\setminus\bigcup_{i=1}^\infty A_i$ is a homeomorphism.
\item $pr^{-1}(A_i)$ is homeomorphic to $Z_i$ for each $i$ (then  for every point $a\in A_i$, the pre-image $pr^{-1}(a)$ must be
homeomorphic to a component of $Z_i$).\qed
\end{lister}
\end{theo}

When $A_1=\{a\}$ is only one ($A_2,A_3,\ldots{}$ are empty) and $Z_1=Z$ is a continuum, we can say the continuum $L(X,Z,a)$ is built
by replacing the point $a\in X$ with~$Z$. The atomic map $pr$ compresses the terminal subcontinuum $Z$ of $L(X,Z,a)$ into~$a$.

\section{Two strategies of proof}

Every diligent schoolboy may prove the following claim.

\begin{prop}
For every zero-dimensional complete space $G$, there is an inverse sequence $(X_n,f^{n+1}_n)_{n=1}^\infty$ such that each space $X_n$
is countable (finite or infinite) and discrete, each bonding map $f^{n+1}_n\colon X_{n+1}\to X_n$ is a surjection, and the inverse
limit\/ $\lim_{\gets} (X_n,f^{n+1}_n)_{n=1}^\infty$ is homeomorphic to~$G$.\qed
\end{prop}

Let $D_0$ be the union of two homeomorphic copies of the bucket-handle continuum, joint at their common end point; $D_0$~is chainable
and has no end point. Let $D_1$ denote a singleton; the point is trivially an end point of~$D_1$. Adikari and Lewis \cite[Theorem
8]{adl-19} have proved that {\em every non-degenerate, hereditarily decomposable chainable continuum has at least one pair of opposite
end points.} It is a good heuristic exercise to construct, for each natural number $k\geq 2$, a Suslinian chainable continuum $D_k$
with exactly $k$ end points (try to modify the $\sin (1/x)$-curve). Adikari and Lewis \cite[pp.\ 621--622]{adl-19} have also
constructed Suslinian chainable continua $D_\infty$ and $D_C$ whose sets of end points are homeo\-morph\-ic, respectively, to the
countably infinite discrete space and the Cantor set. In Remark (A) below we shall see more explicitly what $D_k$ and $D_\infty$ look
like.

\begin{proof}[Proof of Theorem \ref{main}]
For a complete zero-dimensional space $G$, let $(X_n,f^{n+1}_n)_{n=1}^\infty$ be the inverse sequence of Proposition~2. Write $|X_n|$
for the number of elements in $X_n$ or $\infty$ if $X_n$ is infinite. Take the continuum $Y_1=D_{|X_1|}$, and let $h_1\colon X_1\to
Y_1$ be a bijection onto the end point set $E_1$ of~$Y_1$.

\label{prooo} Assume that, for $n=1,\ldots,m$, we have defined a Suslinian chainable continuum $Y_n$ with a one-to-one function
$h_n\colon X_n\to Y_n$ onto the end point set $E_n$ of $Y_n$. Consider the composition $\varphi=h_m f_m^{m+1}\colon X_{m+1}\to E_m$
and, for each $a\in E_m$, put $Z_a=D_{|\varphi^{-1}(a)|}$. We are ready to apply Theorem~4. Let $Y_{m+1}$ be $L(Y_m,Z_a,a)$, and let
$g^{m+1}_m= pr\colon Y_{m+1}\to Y_m$ be the atomic map of Theorem~4; for each $a\in E_m$, $pr^{-1}(a)$ is homeomorphic to
$D_{|\varphi^{-1}(a)|}$. By Proposition 1(d, e) and Ma\'{c}\-ko\-wiak's Corollary, $Y_{m+1}$ is Suslinian, chainable, and its end
point set $E_{m+1}$ consists of all the end points of $pr^{-1}(a)$, $a\in E_m$. Thus, there is a bijection $h_{m+1}\colon X_{m+1}\to
E_{m+1}$ such that $h_{m+1}(\varphi^{-1}(a))= E_{m+1}\cap pr^{-1}(a)$ for each $a\in E_m$. Therefore $g_m^{m+1}h_{m+1}=h_m f_m^{m+1}$.

We have defined an inverse sequence $(Y_n,g^{n+1}_n)_{n=1}^\infty$, and the one-to-one functions $h_n\colon X_n\to Y_n$ form a map
from the inverse sequence $(X_n,f^{n+1}_n)_{n=1}^\infty$ to $(Y_n,g^{n+1}_n)_{n=1}^\infty$; see \cite[pp.\ 101--104]{eng-89} or
\cite[Exercise 2.22]{nad-92}. Take the limit $D_G=\lim_{\gets} (Y_n,g^{n+1}_n)_{n=1}^\infty$ and the induced homeo\-morph\-ic\-al
embedding $h\colon\lim_{\gets} (X_n,f^{n+1}_n)_{n=1}^\infty\to D_G$. The following facts are well-known or easily checked. (1)~{\em An
inverse limit of chainable continua is chainable, too.} (2)~If bonding maps $g^{n+1}_n$ are monotone/atomic, then so are projections
$q_n\colon D_G \to Y_n$. This and Proposition 1(a) yield two more facts. (3)~Since $Y_n$ are Suslinian, $D_G$ is Suslinian, too.
(4)~The set $$E={\bigcap}_{n=1}^\infty q_n^{-1}(E_n)={\lim}_{\gets} (E_n,g^{n+1}_n|E_{n+1})_{n=1}^\infty=h({\lim}_{\gets}
(X_n,f^{n+1}_n)_{n=1}^\infty)$$ consists of end points of $D_G$.\label{proo}
\end{proof}

We shall take a closer look at the simplest continua $D_k$, the continua $Y_n$ of the foregoing proof, and the location of their end
points. The continua will turn out to have a nice three-storey structure, and this will enable us to improve the resulting
continuum~$D_G$.

\begin{rems*}\em
{\bf (A)} The continuum $D_3$, the $\sin (1/x)$-curve, is not arcwise connected, and the complement of its end point set consists of
two arc-components, which are homeomorphic to $(0,1)$. Let $j^3_2\colon D_3\to D_2=[0,1]$ be the atomic map which sends the only
terminal subarc of $D_3$ to~$0$. Let $e_2=0$, $e_2'=1$, and write $e_3,e_3'$ for the ends of $(j^3_2)^{-1}(e_2)$. We obtain $D_{k+1}$
by replacing an end point of $D_k$ with an arc. Thus, there is an atomic map $j^{k+1}_k\colon D_{k+1}\to D_k$ which compresses a
terminal subarc of $D_{k+1}$, sends it to an end point $e_k$ of~$D_k$, and $e_{k+1},e'_{k+1}\in D_{k+1}$ are the ends of the arc
$(j^{k+1}_k)^{-1}(e_k)$. Write
$$D_\omega={\lim}_{\gets} (D_k,j^{k+1}_k)_{k=2}^\infty,\,\, e_\omega=(e_k)_{k=2}^\infty\in D_\omega,\;\mbox{and}\,\, r_k\colon
D_\omega\to D_k$$ for the projections of the inverse limit. $r_k$ are atomic maps for $k=2,3,\ldots{}$ The end point set of $D_\omega$
consists of the points $r_k^{-1}(e_k')$, $k=2,3,\ldots,$ that converge to the end point~$e_\omega$. All the continua lie on the plane.
Finally, we obtain $D_\infty$ when we replace $e_\omega\in D_\omega$ with an arc, i.e.\ there is an atomic map $j^\infty_\omega\colon
D_\infty\to D_\omega$ that compresses a terminal subarc of $D_\infty$ and sends it to $e_\omega$. We can ensure that the sequence of
end points $(r_k j^\infty_\omega)^{-1}(e_k')$ converge to a non-end point of the arc $(j^\infty_\omega)^{-1}(e_\omega)$. Therefore,
the closure of the end point set of $D_\infty$ can be a sequence with its limit---this (the zero-dimensional closure of the end point
set) is the first storey of the structure. What will be important is that, {\em for every atomic map above and below, the pre-image of
any non-end point is a singleton.}

For each of $D_k$, $D_\omega$, $D_\infty$, the complement of the end point set consists of pairwise disjoint arc-components which are
arcs without ends. Write $\mathcal{I}(D_k)$, $\mathcal{I}(D_\omega)$, $\mathcal{I}(D_\infty)$ for the arc-component families.

The third structural element is the following somewhat idiosyncratic property for $X=D_k,D_\omega,D_\infty$:
\begin{lister}
\item[$(\dag)$\hfill]\em For every non-degenerate subcontinuum $P$ of $X$, there is an endless arc $I\in\mathcal{I}(X)$ such that $I\cap P$ is connected and dense in $P$.
\end{lister}
To prove this, see at first that $[0,1]$ with $\mathcal{I}([0,1])=\{(0,1)\}$ has property $(\dag)$. Note that
$\mathcal{I}(\{s\})=\emptyset$ for a singleton $\{s\}$ as $s$ is its end point.

The reader may treat the following as a lemma.

\begin{plainlist*}
\item Claim. \em Since\/ $j^{k+1}_k\colon D_{k+1}\to D_k$ is an atomic map, every point-inverse\/ $(j^{k+1}_k)^{-1}(y)$ of a non-end point\/
$y\in D_k$ is degenerate, and\/ $D_k$ and all point-inverses\/ $(j^{k+1}_k)^{-1}(y)$ have property\/ $(\dag)$---this all implies the
equality
$$\mathcal{I} (D_{k+1})=\{(j^{k+1}_k)^{-1}(I)\colon I \in\mathcal{I}(D_k)\}\cup{\bigcup}_{y\in D_k}\mathcal{I}((j^{k+1}_k)^{-1}(y))$$
and property\/ $(\dag)$ for $D_{k+1}$.
\end{plainlist*}

\begin{proof}
Take a non-end point $x\in D_{k+1}$, and write $y=j^{k+1}_k(x)$. By Proposition 1(e) there are two cases: (1)~$x$ is not an end point
of $(j^{k+1}_k)^{-1}(y)$. Then $(j^{k+1}_k)^{-1}(y)$ is non-degenerate, and a certain arc-component $J$ of the complement of the end
point set of $(j^{k+1}_k)^{-1}(y)$ contains~$x$, i.e.\ $x\in J\in \mathcal{I}((j^{k+1}_k)^{-1}(y))$. From Proposition 1(a) and the
fact that an arc does not contain a non-trivial terminal subcontinuum we infer that every arc $A\ni x$ is contained in
$(j^{k+1}_k)^{-1}(y)$. It follows $J\in\mathcal{I}(D_{k+1})$. (2)~$y$~is not an end point of $D_k$. Then $(j^{k+1}_k)^{-1}(y)$ is a
singleton, and $y\in I\in \mathcal{I}(D_k)$ for a certain $I$. Hence, for every arc $A\ni x$, the restriction $j^{k+1}_k|A$ is an
embedding, and $J=(j^{k+1}_k)^{-1}(I)\in\mathcal{I} (D_{k+1})$.

To prove that $D_{k+1}$ has property $(\dag)$, take a non-degenerate subcontinuum $P$ of $D_{k+1}$. By Proposition 1(a) there are two
cases: (1)~If $P\subset (j^{k+1}_k)^{-1}(y)$, then there is an arc-component $I\in\mathcal{I}((j^{k+1}_k)^{-1}(y))$ which witnesses
property $(\dag)$. (2)~If $P=(j^{k+1}_k)^{-1}(j^{k+1}_k (P))$ and this image of $P$ is non-degenerate, then there is an arc-component
$I\in\mathcal{I}(D_k)$ such that $I\cap j^{k+1}_k (P)$ is connected and dense in $j^{k+1}_k (P)$. Then
$J=(j^{k+1}_k)^{-1}(I)\in\mathcal{I} (D_{k+1})$ is homeomorphic with $I$, $j^{k+1}_k$ homeo\-morph\-ic\-al\-ly maps $J\cap P$ onto
$I\cap j^{k+1}_k (P)$, and hence $J\cap P$ is connected. The set $\cl (J\cap P)$ is a continuum with $j^{k+1}_k (\cl (J\cap P))=\cl
j^{k+1}_k (J\cap P)=j^{k+1}_k (P)$, and, as $j^{k+1}_k$ is atomic, we obtain $\cl (J\cap P)=P$. Therefore, $D_{k+1}$ has property
$(\dag)$.
\end{proof}

We have shown by induction that $D_k$, $k=2,3,\ldots,$ satisfy~$(\dag)$. Now, we shall prove $(\dag)$ for an inverse limit.

\begin{plainlist*}
\item Claim. \em Since, for $k=2,3,\ldots,$ continua $D_k$ and maps $j^{k+1}_k\colon D_{k+1}\to D_k$ satisfy the assumptions of the previous Claim, we
obtain for $D_\omega={\lim}_{\gets} (D_k,j^{k+1}_k)_{k=2}^\infty$ the equality
$$\mathcal{I} (D_\omega)={\bigcup}_{k=2}^\infty\{r_k^{-1}(I)\colon I\in\mathcal{I}(D_k)\}$$ and property $(\dag)$.
\end{plainlist*}

\begin{proof}
Using the statement (b) of Theorem~2, we can easily check for a point $x=(x_k)_{k=2}^\infty\in D_\omega$ that, if $x_k$ is an end
point of $D_k$ for $k=2,3,\ldots{}$, then $x$ is an end point of $D_\omega$. Therefore, if $x$ is not an end point of $D_\omega$, then
there is a $k$ such that $x_k,x_{k+1},\ldots$ are not end points of $D_k,D_{k+1},\ldots{}$, respectively. Then $x_k\in
I\in\mathcal{I}(D_k)$ and, by the same argument as in the previous proof applied to $r_k\colon D_\omega\to D_k$, we obtain that $x\in
J=r_k^{-1}(I)\in\mathcal{I}(D_\omega)$.

To prove that $D_\omega$ satisfies $(\dag)$, take a non-degenerate subcontinuum $P$ of $D_\omega$ and a $k$\/ such that $r_k(P)$ is
non-degenerate. There is an endless arc $I\in\mathcal{I}(D_k)$ such that $I\cap r_k(P)$ is connected and $r_k(P)\subset\cl (I\cap
r_k(P))$. By the same argument as in the previous proof, $J=r_k^{-1}(I)\in\mathcal{I}(D_\omega)$ witnesses property $(\dag)$ for
$D_\omega$.
\end{proof}

In the same way as above, since $j^\infty_\omega$ is an atomic map, {\em the complement of the end point set of $D_\infty$ is a
disjoint union of arc-components in $\mathcal{I}(D_\infty)$ homeomorphic to $(0,1)$, and $D_\infty$ has property\/ $(\dag)$.}

{\bf (B)} Let us come back to the proof of Theorem 1 and observe that we can make the continua $Y_n$ and its inverse limit so that
they have the three-storey structure described above. Let us take the simplest $D_\infty$, whose end points form a sequence convergent
to a non-end point. Also note that each restriction $g^{n+1}_n|(g^{n+1}_n)^{-1}(Y_n\setminus E_n)\colon (g^{n+1}_n)^{-1}(Y_n\setminus
E_n)\to Y_n\setminus E_n$ is a homeo\-morph\-ism. It is clear by induction that the accumulation points of $\cl E_n$ are not end
points of~$Y_n$. Hence, $\cl E_n$ is countable for each~$n$. It follows that
$$\cl E={\bigcap}_{n=1}^\infty q_n^{-1}(\cl E_n)={\lim}_{\gets} (\cl E_n,g^{n+1}_n|\cl E_{n+1})_{n=1}^\infty$$ is zero-dimensional.
Secondly, the arc-components of $Y_n\setminus E_n$ and $D_G\setminus E$ are homeo\-morph\-ic to $(0,1)$ by induction and the same
argument as in Remark (A), and they form countable families $\mathcal{I}(Y_n)$ and $\mathcal{I}(D_G)$. Therefore, {\em $D_G$ is the
union of a closed zero-dimensional set\/ $\cl E$ and the countable family $\mathcal{I}(D_G)$ of pairwise disjoint subsets homeomorphic
to $(0,1)$.} Finally, by induction and the same argument as in Remark~(A) {\em each $Y_n$ and $D_G$ have property\/ $(\dag)$.}

{\bf (C)} Having noted the structure of the continuum $D_G$, we can repeat Ja\-ni\-szew\-ski's procedure \cite{jan-13} exactly without
changes: into each endless arc of $D_G$ we insert an arc as a terminal subcontinuum, in this way we damage all endless arcs of $D_G$,
and then, proceeding to infinity, we damage all newly appearing arcs in order to obtain {\em an arcless, Suslinian, chainable
continuum $D_G'$ whose end points form a set homeo\-morph\-ic to~$G$.} This can readily be formalised as an inverse limit.

{\bf (D)} A [{\bf weak}] {\bf Cook} continuum is one, call it $X$, such that each non-constant map from a subcontinuum $P$ of $X$ to
$X$ is the identity $\id_P\colon P\to X$ [respectively: has $P\cap f(P)\neq\emptyset$]. Ma\'{c}kowiak's aim \cite{mac-85} was to
construct a Suslinian, chainable weak Cook cotinuum (his terminology was somewhat different). Let us consider the continuum $D_G$ of
Remark (B), and let $C=\{c_i\}_{i=1}^\infty\subset D_G\setminus E$ be a set such that, for each endless arc $I\in\mathcal{I} (D_G)$,
the intersection $I\cap C$ is dense in~$I$. Let $Z_i$, $i=1,2,\ldots,$ be pairwise disjoint non-degenerate subcontinua of
Ma\'{c}kowiak's weak Cook continuum. Finally, let $D_G''$ be the continuum $L(D_G,Z_i,c_i)$ of Theorem 4 and let $pr\colon D_G''\to
D_G$ be that atomic map. By Proposition 1(d, e) and the Corollary to Theorem~3, $D_G''$ is chainable, Suslinian, and its end point set
is $pr^{-1}(E)$, which is homeomorphic to~$G$.

\begin{plainlist*}
\item Claim. \em $D_G''$ is a weak Cook continuum.
\end{plainlist*}

\begin{proof}
Let us start with showing that {\em every map $f\colon P\to D_G''$ from a subcontinuum $P$ of $D_G$ is constant.} If $P$~is
non-degenerate, then by property $(\dag)$ there is an endless arc $I\in\mathcal{I}(D_G)$ such that $I\cap P$ is connected and
$P\subset\cl (I\cap P)$. Let $A\subset I\cap P$ be an arc with ends. There are two cases. (1) $Q=f(A)\subset pr^{-1}(y)$ for a point
$y\in D_G$. If $f(A)$ were non-degenerate, then $y$ would be a $c_i$ and $Z_i$ would contain an arc, which is impossible for a weak
Cook continuum. Thus, $f|A$ is constant. (2)~The image $pr(Q)$ is non-degenerate and $Q =pr^{-1}(pr(Q))$ by Proposition 1(a). It
follows from property $(\dag)$ applied to $pr(Q)$ that a certain $c_i\in pr(Q)$. Thus, $pr^{-1}(c_i)\subset Q$. Theorem 12.46 in
\cite{nad-92} says that {\em each map onto a chainable continuum is weakly confluent\/\footnote{A map $f\colon X\to Y$ is {\bf weakly
confluent} if, for every subcontinuum $Q$ of $Y$, there is a subcontinuum $P$ of $X$ with $f(P)=Q$.}}, i.e.\ there is a continuum
$R\subset A$ such that $f(R)= pr^{-1}(c_i)$. Hence, again $Z_i$ contains an arc. A contradiction. Therefore, the case (2) does not
hold, and $f(A)$ is degenerate. It follows that $f|(I\cap P)$ is constant, and $f$~is constant as $I\cap P$ is dense in~$P$.

Now, take a subcontinuum $P$ of $D_G''$ and a map $f\colon P\to D_G''$ with $P\cap f(P)=\emptyset$. We shall repeat the above schema
to prove that, {\em if $pr(P)$ is a singleton, then $f$~is a constant map.} Assume $pr|P$ is constant. There are two cases. (1)~If
$Q=f(P)$ is contained in a point-inverse $pr^{-1}(y)$, $y\in D_G$, then either $f$ is constant or $f(P)\subset pr^{-1}(c_i)$. In case
$f(P)\subset pr^{-1}(c_i)$, $f$ induces a map into the subcontinuum $Z_i$ of Ma\'{c}kowiak's continuum, and hence $f$ is constant.
(2)~The image $pr(Q)$ is non-degenerate and $pr^{-1}(pr(Q))=Q$. Since $D_G$ has property $(\dag)$, there is an endless arc
$I\in\mathcal{I}(D_G)$ such that $I\cap pr(Q)$ is dense in $pr(Q)$. Thus, there is a  $c_i\in I\cap pr(Q)$ and $pr^{-1}(c_i)\subset
Q=f(P)$. By \cite[Theorem 12.46]{nad-92}, $f$~is weakly confluent, i.e.\ there is a continuum $R\subset P$ such that $f(R)=
pr^{-1}(c_i)$. As $R$ and $pr^{-1}(c_i)$ are disjoint subcontinua of Ma\'{c}kowiak's continuum, $f|R$ is constant, and the case (2)
does not hold, in fact.

To continue assume $P=pr^{-1}(pr(P))$. Then, for every $y\in pr(P)$, the restriction $f|pr^{-1}(y)$ is constant. It follows that $f$
has a factorisation $f=g\,pr$, where $g\colon pr(P)\to D_G''$ is defined by the correct formula $g(y)=f(pr^{-1}(y))$. By the first
paragraph of this proof, $g$ is constant, and so is~$f$.
\end{proof}

Therefore, {\em $D_G''$ is a weak Cook, Suslinian chainable continuum whose set of end points is homeomorphic to~$G$.}

{\bf (E)} In \cite[Footnote 6]{krz-10} we observed that {\em no non-degenerate Suslinian continuum is a Cook continuum.}
\end{rems*}

In our proof of Theorem~1 we have multiplied end points in~$E_n$ in order to breed the $E_n$s' inverse limit. There is another
strategy, although it will not be entirely effective. We can take Adikari and Lewis' continuum $D_C$ whose end points form a Cantor
set, and destroy those unwanted. Again the next proof is left to the reader.

\begin{prop}
Any\/ $\sigma$-compact zero-dimensional space is the union of a sequence of pairwise disjoint compact subsets $F_1,\ldots,F_i,\ldots$
with\/ $\lim_{i\to\infty}\diam F_i=0$.\qed
\end{prop}

\begin{proof}[Unsuccessful attempt of proof of Theorem \ref{main}]
The end points of the chainable continuum $D_C$ are in the Cantor set $C\subset D_C$. Given a complete zero-dimensional space $G$, we
can assume that $G$ is a $G_\delta$ subset of~$C$. Thus, $C\setminus G$ is the union $\bigcup_{i=1}^\infty F_i$ of Proposition~3. We
use $D_0$, a chainable continuum without end points, and we take Pol and Re\'{n}ska's continuum $D_G=L(D_C,F_i\times D_0,F_i)$ with
the atomic map $pr\colon D_G\to D_C$. Since $D_0$ has no end point, $pr^{-1}(G)$ is the end point set of $D_G$ by Proposition 1(e).
The restriction $pr|pr^{-1}(G)\colon pr^{-1}(G)\to G$ is a homeomorphism. $D_G$ is chainable by Ma\'{c}\-ko\-wiak's Corollary.
Unfortunately, $D_G$ contains homeomorphic copies of the indecomposable bucket-handle since $D_0$ does.
\end{proof}

\section{Questions and reflections}

Let us restate Adikari and Lewis' problem. (Note that every point of a pseudo-arc is its end point.)

\begin{prob}
Construct chainable continua with interesting end point sets or disprove their existence.

In particular, does there exist a hereditarily decomposable {\em [}respectively: Suslin\-ian\/{\em ]} chainable continuum whose set of
end points is one-dimensional? Can Kuratow\-ski's set\/\footnote{See \cite[Problem 1.2.E]{eng-95}.} or a complete Erd\H{o}s
space\/\footnote{See J.\ J.\ Dijkstra and J.\ van Mill\/ \cite{dvm-09}, \cite[Problem 6.3.24]{eng-89}, and  cf.\ \cite[Example
1.2.15]{eng-95}.} be the set of end points of a hereditarily decomposable {\em [}respectively: Suslinian\/{\em ]} chainable continuum?
\end{prob}

If we take the Knaster-Kuratowski fan\footnote{See \cite[Problem 6.3.23]{eng-89} and \cite[Problem 1.4.C]{eng-95}.} and remove its
dispersion point, then we obtain a one-dimensional hereditarily disconnected space $M$ that cannot be contained in the end point set
of an hereditarily decomposable chainable continuum. Indeed, E.\ Pol \cite{pol-78} proved that {\em every completely metrisable space
which contains~$M$ also contains an arc,}\/ and thus the end point set would contain an arc.

Ma\'{c}kowiak's Corollary to Theorem~3 may be considered as the inverse invariance of chainability under atomic maps. But it is not
entirely general. Many authors (the present one in this number) construct their examples by inserting a chainable continuum into
another chainable continuum instead of its point. It is intriguing that we do not know whether we always obtain a chainable continuum.
When we insert the continuum into a hereditarily decomposable continuum, that Corollary says `Yes'! In a joint paper \cite{chk-10}
with M.\ G.\ Charalambous we have proved that, when we insert a pseudo-arc into a pseudo-arc, we again obtain a pseudo-arc. Therefore,
we consider the following problem as most important.

\begin{prob}
Suppose that $X$ is a continuum, it contains a terminal chainable subcontinuum $P$, and the quotient continuum $X/P$ is chainable.
Under what circumstances is~$X$ itself chainable? If moreover $X/P$ is the bucket-handle continuum?
\end{prob}

We have to assume that the quotient map is atomic as the simple example of a monotonic map from the figure T continuum to $[0,1]$
shows. Much less obvious example of a non-terminal subcontinuum is that of J.\ F.\ Davis and W.\ T.\ Ingram \cite{dai-88}; in that
case the quotient continuum is indeed the bucket-handle. In light of Fugate's theorem it is sufficient to solve the problem for
indecomposable quotients $X/P$.

\begin{center}
*\end{center}

The author would never expect that property $(\dag)$ might be advantageous. I~found it when I tried to prove that the continuum
$D_G''$ of Remark (D) is weak Cook. That was analysis, the method to find cause or to describe conducive circumstances when we know
the result. Then, that was only formalism which led me to find appropriate hypotheses for the claims on atomic maps of Remark (A) and
simpler and simpler proofs. These are not deep remarks, but we should illustrate to our students that analysis, not formal logic, is a
method of the {\em ars inveniendi.} Synthesis is a method of the {\em ars iudicandi.} Between them there is the art of establishing
convenient notation and terminology, and the three: heuristic, logic, and semantics, play complementary and supplementary roles. This
is, in fact, the old division of the trivium: rhetoric, dialectic, and grammar. And it is that disfavoured rhetoric that is pertinent
to the quest for the shape of the universe or for the nature of the human brain's functioning. {\em Matematyka prowadzi nas ku
\'{s}wiatu}\/ (Polish) [Mathematics guides us towards the world]---this is a dictum of Professor Mioduszewski.


\begin{thebibliography}{33}

\bibitem{adl-19}\label{adl}
R.\ Adikari and W.\ Lewis, {\em Endpoints of nondegenerate hereditarily decomposable chainable continua,} Houston J.\ Math.\ 45
(2019), no.\ 2, 609--624.

\bibitem{ach-59}
R.\ D.\ Anderson, G.\ Choquet, {\em A plane continuum no two of whose nondegenerate subcontinua are homeomorphic: an application of
inverse limits,} Proc.\ Amer.\ Math.\ Soc.\ 10 (1959), no.\ 3, 347--353.

\bibitem{bing-51}
R.\ H.\ Bing, {\em Snake-like continua,} Duke Math.\ J.\ 18 (1951), no.\ 3, 653--663.

\bibitem{chk-10}
M.\ G.\ Charalambous and J.\ Krzempek, {\em On Dimensionsgrad, resolutions, and chainable continua,} Fund.\ Math.\ 209 (2010), no.\ 3,
243--265.

\bibitem{dai-88}
J.\ F.\ Davis and W.\ T.\ Ingram, {\em An atriodic tree-like continuum with positive span which admits a monotone mapping to a
chainable continuum,} Fund.\ Math.\ 131 (1988), no.\ 1, 13--24.

\bibitem{dvm-09}
J.\ J.\ Dijkstra and J.\ van Mill, {\em Characterizing complete Erd\H{o}s space,} Canad.\ J.\ Math.\ 61, 2009, no.\ 1, 124--140.

\bibitem{fug-69}
J.\ B.\ Fugate, {\em A characterization of chainable continua,} Canadian J.\ Math.\ 21 (1969), 383--393.

\bibitem{eng-89}
R.\ Engelking, {\em General Topology,} Heldermann Verlag, Berlin, 1989.


\bibitem{eng-95}
R.\ Engelking, {\em Theory of Dimensions, Finite and Infinite,} Heldermann Verlag, Lemgo, 1995.

\bibitem{isb-59}
J.\ R.\ Isbell, {\em Embeddings of inverse limits,} Ann.\ of Math.\ (2) 70 (1959), no.\ 1, 73--84.

\bibitem{jan-13}
Z.\ Janiszewski, {\em \"{U}ber die Begriffe ``Linie'' und ``Fl\"{a}che'',} in: {\em Proceedings of the Fifth International Congress of
Mathematicians (Cambridge, 22--28 August 1912),} Vol.\ II, 126--128, Cambridge University Press, Cambridge, 2013; also in: Z.\
Ja\-ni\-szew\-s\-ki, {\em Oeuvres choisies,} 127--129, PWN, Warszawa, 1962.

\bibitem{krz-10}
J.\ Krzempek, {\em Fully closed maps and non-metrizable higher-dimensional Anderson-Choquet continua,} Colloq.\ Math.\ 120 (2010),
no.\ 2, 201--222.

\bibitem{mac-85}
T.\ Maćkowiak, {\em The condensation of singularities in arc-like continua,} Houston J.\ Math.\ 11 (1985), no.\ 4, 535--558.

\bibitem{mio-62}
J.\ Mioduszewski, {\em A functional conception of snake-like continua,} Fund.\ Math.\ 51 (1962/63), no.\ 2, 179--189.

\bibitem{mio-64}
J.\ Mioduszewski, {\em Everywhere oscillating functions, extension of the uniformization and homogeneity of the pseudo-arc,} Fund.\
Math.\ 56 (1964), no.\ 2, 131--155.

\bibitem{nad-92}
S.\ B.\ Nadler, Jr.\ {\em Continuum Theory. An Introduction,} Monographs and Textbooks in Pure and Applied Mathematics 158, Marcel
Dekker, Inc., New York, 1992.

\bibitem{pol-78}
E.\ Pol, {\em Strongly metrizable spaces of large dimension each separable subspace of which is zero-dimensional,} Colloq.\ Math.\ 39
(1978), no.\ 1, 25--27, 189.

\bibitem{por-03}
E.\ Pol, M.\ Reńska, {\em On Bing points in infinite-dimensional hereditarily indecomposable continua,} Topology Appl.\ 123 (2003),
no.\ 3, 507--522.

\bibitem{why-30}
G.\ T.\ Whyburn, {\em A continuum every subcontinuum of which separates the plane,} Amer.\ J.\ Math.\ 52 (1930), no.\ 2, 319--330.

\end{thebibliography}
\end{document}